\documentclass{article}

\usepackage[compress]{cite}
\usepackage{amssymb}
\usepackage{amsthm}
\usepackage{dsfont}
\usepackage{mathtools}
\usepackage{amsmath}
\usepackage{xfrac}
\usepackage{stmaryrd}
\usepackage[hyphenbreaks]{breakurl}
\usepackage[hyphens]{url}
\usepackage{tikz}
\usepackage[margin=1in]{geometry}

\usepackage[font=footnotesize,labelfont=it]{caption}

\renewcommand{\v}[1]{\mathbf{#1}}

\newtheorem{definition}{Definition}

\newtheorem{remark}{Remark}

\newtheorem{assumption}{Assumption}
\newtheorem{theorem}{Theorem}

\newcommand{\revise}[1]{{\color{black}#1}}

\DeclareMathOperator*{\argmin}{arg\,min}

\begin{document}
	
	\title{Deconvoluting Kernel Density Estimation and Regression for Locally Differentially Private Data} 
	
	\author{Farhad Farokhi 
		\thanks{The author is with the Department of Electrical and Electronic Engineering at the University of Melbourne. e-mail: farhad.farokhi@unimelb.edu.au }
	}

\date{}

\maketitle

\begin{abstract} Local differential privacy has become the gold-standard of privacy literature for gathering or releasing sensitive individual data points in a privacy-preserving manner. However, locally differential data can twist the probability density of the data because of the additive noise used to ensure privacy. In fact, the density of privacy-preserving data (no matter how many samples we gather) is always flatter in comparison with the density function of the original data points due to convolution with privacy-preserving noise density function. The effect is especially more pronounced when using slow-decaying privacy-preserving noises, such as the Laplace noise. This can result in under/over-estimation of the heavy-hitters. This is an important challenge facing social scientists due to the use of differential privacy in the 2020 Census in the United States. In this paper, we develop density estimation methods  using smoothing kernels. We use the framework of deconvoluting kernel density estimators  to remove the effect of privacy-preserving noise. This approach also allows us to adapt the results from non-parametric regression with errors-in-variables to develop regression models based on locally differentially private data.  We demonstrate the performance of the developed methods on financial and demographic datasets. 
\end{abstract}

\section*{Introduction}
Government regulations, such as the roll-out of the General Data Protection Regulation in the European Union (EU)\footnote{\url{https://gdpr-info.eu}}, the California Consumer Privacy Act\footnote{\url{https://oag.ca.gov/privacy/ccpa}}, and the development of the Data Sharing and Release Bill in Australia\footnote{\url{https://www.pmc.gov.au/public-data/data-sharing-and-release-reforms}} increasingly prohibit sharing customer’s data without explicit consent~\cite{bennett2020revisiting}. 

A strong candidate for ensuring privacy is differential privacy. Differential privacy intuitively uses randomization to provide plausible deniability for the data of an individual by ensuring that the statistics of privacy-preserving outputs do not change significantly by varying the data of an individual~\cite{Dwork2008,dwork2006calibrating}. Companies like Apple\footnote{\url{https://www.apple.com/privacy/docs/Differential_Privacy_Overview.pdf}}, Google\footnote{\url{https://developers.googleblog.com/2019/09/enabling-developers-and-organizations.html}},  Microsoft\footnote{\url{https://www.microsoft.com/en-us/ai/ai-lab-differential-privacy}}, and LinkedIn\footnote{\url{https://engineering.linkedin.com/blog/2019/04/privacy-preserving-analytics-and-reporting-at-linkedin}} have rushed to develop projects and to integrate differential privacy into their products. Even, the US Census Bureau has decided to implement differential privacy in 2020 Census~\cite{abowd2018us}. Of course, this has created much controversy pointing to ``ripple effect on the many public and private organizations that conduct surveys based on census data"~\cite{Mervis114}. 

A variant of differential privacy is local differential privacy in which all data points are randomized before being used by the aggregator, who attempts
to infer the data distribution or some of its properties~\cite{dewri2013local,6686179,kairouz2014extremal}. This is in contrast with differential privacy in which the data is first processed and then obfuscated by noise. Local differential privacy ensures that the data is kept private from the aggregator by adding noise to the individual data entries before the aggregation process. This is a preferred choice when dealing with untrusted aggregators, e.g., third party service providers or  commercial retailers with financial interests, or when it is desired to release an entire dataset publicly for research in a privacy-preserving manner~\cite{liu2020local}. \revise{Differential privacy is in spirit close to randomized response methods introduced originally in~\cite{warner1965randomized} to reduce potential bias due to non-response and social desirability when asking questions about sensitive topics. The randomized response can be used to conceal individual responses (i.e., protect individual privacy) so that the respondents are more inclined to answer truthfully~\cite{kuk1990asking,blair2015design,boruch1971assuring,fox1984measuring}. In fact, for questions with binary answer, the randomized response method with forced response (i.e., a respondent determines whether to  answer a sensitive question truthfully or with forced yes/no based on flipping a biased coin) is differentially private and the probability of the head for the coin determines the privacy budget in differential privacy~\cite{kearns2019ethical}. However, differential is a more general and flexible methodology that can be used for categorical and non-categorical (i.e., continuous domain) questions~\cite{abadi2016deep,friedman2010data,abowd2018us}. This paper specifically consider the problem of analyzing privacy-preserving data on continuous domain, which is out of the scope of randomized response methodology. }

Locally differential data can significantly distort our estimates of the probability density of the data because of the additive noise used to ensure privacy. The density of privacy-preserving data can become flatter in comparison with the density function of the original data points due to convolution of its density with privacy-preserving noise density. The situation can be even more troubling when using slow-decaying privacy-preserving noises, such as the Laplace noise. This concern is true irrespective of how many samples are gathered. This can result in under/over-estimation of the heavy-hitters, a common and worrying criticism of using differential privacy in the US Census~\cite{NYtimesCensus}. 

Estimating probability distributions/densities under differential privacy is of extreme importance as it is often the first step in gaining more important insights into the data, such as regression analysis. However, most of the existing work on probability distributions estimation based on locally differential private data focuses on 
categorical data~\cite{acharya2019hadamard,bassily2015local,erlingsson2014rappor, wang2017locally,ye2018optimal}. For categorical data (in contrast with numerical data), the privacy-preserving noise is no longer additive, e.g., the so-called exponential mechanism~\cite{mcsherry2007mechanism} or other boutique differential privacy mechanisms~\cite{li2020estimating} are often employed that are not on the offer in the 2020 US Census. \revise{The density estimation results for categorical data are also related to de-noising results in randomized response methods~\cite{blair2015design}.} The work on continuous domains is often done by binning or quantizing the domain. However, finding the optimal number of bins or quantization resolution depending on privacy parameters, data distribution, and number of data points is a challenging task.

In this paper, we take a different approach to density estimation by using kernels and thus eliminating the need to quantize the domain. \revise{Kernel density estimation is a non-parametric way to estimate the probability density function of a random variable using its samples proposed independently by Parzen~\cite{Parzen1962} and Rosenblatt~\cite{Rosenblatt1956}. This methodology was extended to multi-variate variables in~\cite{murthy1966nonparametric,loftsgaarden1965}. These estimators work based on batches of data; however, they can also be made recursive~\cite{mokkadem2009stochastic,slaoui2018bias, slaoui2019data}. When the data samples are noisy because of measurement noise or, as in the case of this paper, privacy-preserving noise, we need to eliminate the effect of the additive noise kernel density estimation by deconvolution~\cite{carroll1988optimal}. Therefore,} we use the framework of deconvoluting kernel density estimators~\cite{stefanski1990deconvolving,neumann1997effect, delaigle2008deconvolution,carroll1988optimal} to remove the effect of privacy-preserving noise, which is often in the form of Laplace noise~\cite{dwork2014algorithmic}. This approach also allows us to adapt the results from non-parameteric regression with errors-in-variables~\cite{delaigle2007nonparametric,ioannides1997nonparametric,fan1993nonparametric} to develop regression models based on locally differentially private data. \revise{This is the first time that  deconvoluting kernel density estimators have been used for analyze differentially-private data.} This is an important challenge facing social science researchers and demographers following the changes administered in the 2020 Census in the United States~\cite{abowd2018us}. 

\section*{Methods}
Consider independently distributed data points $\{\v{x}[i]\}_{i=1}^n\subset\mathbb{R}^q$, for some fixed dimension $q\geq 1$, from common probability density function $\phi_{\v{x}}$. Each data point $\v{x}[i]\in\mathbb{R}^q$ belongs to an individual. Under no privacy restrictions, the data points can be provided to the central aggregator to construct an estimate of the density $\phi_{\v{x}}$ denoted by $\widehat{\phi}_{\v{x}}$. We may use kernel $K$, which is a bounded even probability density function, to generate the density estimate $\widehat{\phi}_{\v{x}}$. A widely recognized example of a kernel is the Gaussian kernel~\cite{wand1994kernel} in
\begin{align} \label{eqn:Gaussian_kernel}
	K({\v{x}})=\frac{1}{\sqrt{(2\pi)^q}}\exp\left(-\frac{1}{2}\v{x}^\top \v{x}\right).
\end{align}
In the big data regime $n\gg 1$, the choice of the kernel is not crucial to the accuracy of kernel density estimators so long as it meets the conditions in~\cite{stefanski1990deconvolving}. In this paper, we keep the kernel general. By using kernel $K$, we can construct the estimate
\begin{align} \label{eqn:density_no_privacy}
	\widehat{\phi}^{\mathrm{np}}_{\v{x}}(\v{x})=\frac{1}{nh^q}\sum_{i=1}^{n}K((\v{x}-\v{x}[i])/h),
\end{align}
where $h>0$ is the bandwidth. The bandwidth is often selected such that $h\rightarrow 0$ as $n\rightarrow \infty$. The optimal rate of decay for the bandwidth has been established for families of distributions~\cite{stefanski1990deconvolving,carroll1988optimal}.

\revise{\begin{remark}
The problem formulation in this paper considers real-valued data as opposed as categorical data. This distinguishes the paper from the computer science literature on this topic, which primarily focuses on 
categorical data~\cite{acharya2019hadamard,bassily2015local,erlingsson2014rappor, wang2017locally,ye2018optimal}. Real-valued data can arise in two situations. First, the posed question can be non-categorical, e.g., credit rating for loans or the interest rates of loans. We will consider this in one of our experimental results. However, aggregated categorical data can also be real-valued. For instance, the 2020 US Census reports the aggregate number of individuals from a race or ethnicity group within different counties. These numbers will be made differentially private as part of the US Census Bureau's privacy initiative~\cite{abowd2018us}. Therefore, the methods developed in this paper are still relevant to categorical data, albeit in aggregated forms.
\end{remark}
}

As discussed in the introduction, due to privacy restrictions,  the exact data points $\{\v{x}[i]\}_{i=1}^n$ might not be available to generate the density estimate in~\eqref{eqn:density_no_privacy}. The aggregator may only have access to noisy versions of these data points:
\begin{align} \label{eqn:additive}
	\v{z}[i]=\v{x}[i]+\v{n}[i],
\end{align}
where $\v{n}[i]$ is a privacy-preserving additive noise. To ensure differential privacy, Laplace additive noises is often used~\cite{dwork2014algorithmic}. For any probability density $\phi$, we use the notation $\mathrm{supp}(\phi)$ to denote its support set, i.e., $\mathrm{supp}(\phi):=\{\xi:\phi(\xi)>0\}$. 

\begin{assumption}[Bounded Support] \label{assum:1}
	$\mathrm{supp}(\phi_{\v{x}})\subseteq\prod_{i=1}^q [\underline{x}_i,\overline{x}_i]$ for finite constants $\underline{x}_i\leq \overline{x}_i$. 
\end{assumption}

Assumption~\ref{assum:1} is without loss of generality as we are always dealing with bounded domains in social sciences with \textit{a priori} known bounds on the data (e.g., the population of a region). 

\begin{definition}[Local Differential Privacy] \label{def:LDP}
	The reporting mechanism in~\eqref{eqn:additive} is $\epsilon$-(locally) differentially private for $\epsilon\geq 0$ if 
	\begin{align*}
		\mathbb{P}\{\v{x}[i]+\v{n}[i]\in\mathcal{Z}|\v{x}[i]=\v{x}\}\leq \exp(\epsilon) \mathbb{P}\{\v{x}[i]+\v{n}[i]\in&\mathcal{Z}|\v{x}[i]=\v{x}'\},\quad \forall \v{x},\v{x}'\in\mathrm{supp}(\phi_{\v{x}}),
	\end{align*}
	for any Borel-measurable set $\mathcal{Z}\subseteq \mathbb{R}^q$. 
\end{definition}

Definition~\ref{def:LDP} ensures that the statistics of privacy-preserving output $\v{x}[i]+\v{n}[i]$, determined by its distribution, do not change ``significantly'' (the magnitude of change is bounded by the privacy parameter $\epsilon$) if the data of individual $\v{x}[i]$ changes. If $\epsilon\rightarrow 0$, the output becomes more noisy and a higher privacy guarantee is achieved. Laplace additive noise is generally used to ensure differential privacy. This is formalized in the following theorem, which is borrowed from~\cite{dwork2014algorithmic}.

\begin{theorem} \label{tho:1} Let $\{\v{n}[i]\}_{i=1}^n$ be distributed according to the common multivariate Laplace density:
	\begin{align*}
			\phi_{\v{n}}(\v{n})=\frac{1}{2^q\prod_{j=1}^q b_j}\exp\left(-\sum_{j=1}^q\frac{|n_j|}{b_j}\right),
	\end{align*}
	where $n_j$ is the $j$-th component of $\v{n}\in\mathbb{R}^q$. The reporting mechanism in~\eqref{eqn:additive} is $\epsilon$-locally differentially private if $b_j=q(\overline{x}_j-\underline{x}_j)/\epsilon$ for $j\in\{1,\dots,q\}$. 
\end{theorem}

In what follows, we assume that the reporting policy in Theorem~\ref{tho:1} is used to generate locally differentially private data points.  Since $\{\v{n}[i]\}_{i=1}^n$ are distributed according to the common density $\phi_{\v{n}}(\v{n})$,  $\{\v{z}[i]\}_{i=1}^q$ would also follow a common probability density, which is denoted by $\phi_{\v{z}}$. Note that 
\begin{align} \label{eqn:zvn}
	\Phi_{\v{z}}(\v{t})=\Phi_{\v{x}}(\v{t})\Phi_{\v{n}}(\v{t}),
\end{align}
where $\Phi_{\v{z}}$, $\Phi_{\v{x}}$, and $\Phi_{\v{n}}$ are the characteristic functions of $\phi_{\v{z}}$, $\phi_{\v{x}}$, and $\phi_{\v{n}}$. Using~\eqref{eqn:zvn},  we can use any approximation of $\Phi_{\v{z}}$ to construct an approximation of $\Phi_{\v{x}}$ and thus estimate $\phi_{\v{x}}$. 
If we use kernel $K$ for estimating density of $\v{z}[i]$, $\forall i$, we get
\begin{align*}
	\widehat{\phi}_{\v{z}}(\v{z})=\frac{1}{nh^q}\sum_{i=1}^{n}K((\v{z}-\v{z}[i])/h).
\end{align*}
Here, $\widehat{\phi}_{\v{z}}$ is used to denote the approximation of $\phi_{\v{z}}$. The characteristic function of  $\widehat{\phi}_{\v{z}}$ is given by
\begin{align*}
	\widehat{\Phi}_{\v{z}}(\v{t})
	=&\Phi_{K}(h\v{t})\widehat{\Phi}(\v{t}),
\end{align*}
where $\Phi_{K}(\v{t})$ is the characteristic function of $K$ and $\widehat{\Phi}(\v{t})$ is the empirical characteristic function of measurements $\{\v{z}[i]\}_{i=1}^n$, defined as
\begin{align*}
	\widehat{\Phi}(\v{t})=\frac{1}{n}\sum_{i=1}^{n}\exp\left(i \v{t}^\top \v{z}[i] \right).
\end{align*}
Therefore, the characteristic function of  $\widehat{\phi}_{\v{x}}$ is given by
\begin{align*}
	\widehat{\Phi}_{\v{x}}(\v{t})=\frac{\Phi_{K}(H\v{t})\widehat{\Phi}(\v{t})}{\Phi_{\v{n}}(\v{t})}
\end{align*}
Further, note that
\begin{align*}
	\Phi_{\v{n}}(\v{t})
	=&\mathbb{E}\left\{\exp\left(i \v{t}^\top \v{n}\right)\right\}\\
	=&\mathbb{E}\left\{\exp\left(i t_{1} n_{1}\right)\exp\left(i t_{2} n_{2}\right)\cdots \exp\left(i t_{q} n_{q}\right)\right\}\\
	=&\mathbb{E}\left\{\exp\left(i t_{1} n_{1}\right)\right\}\mathbb{E}\left\{\exp\left(i t_{2} n_{2}\right)\right\}\cdots \mathbb{E}\left\{\exp\left(i t_{q} n_{q}\right)\right\}\\
	=&\prod_{j=1}^q \frac{1}{1+b_j^2 t_j^2},
\end{align*}
where $t_j$ is the $j$-th component of $\v{t}\in\mathbb{R}^q$.  We get
\begin{align} \label{eqn:density_private}
	\widehat{\phi}_{\v{x}}(\v{x})=\frac{1}{nh^q} \sum_{i=1}^n \widehat{K}_h((\v{x}-\v{z}[i])/h),
\end{align}
where
\begin{align*}
	\widehat{K}_h(\v{x})
	&=\frac{1}{(2\pi)^q}\int_{\mathbb{R}^q}\exp(-i \v{t}^\top \v{x})\frac{\Phi_{K}(\v{t})}{\Phi_{\v{n}}(\v{t}/h)}\mathrm{d}\v{t}\\
	&=\frac{1}{(2\pi)^q}\int_{\mathbb{R}^q}\exp(-i \v{t}^\top \v{x})\prod_{j=1}^q \left(1+\frac{b^2}{h^2} t_j^2\right)\Phi_{K}(\v{t})\mathrm{d}\v{t}\\
	&=\prod_{j=1}^q \left(1-\frac{b_j^2}{h^2} \frac{\partial^2}{\partial x_j^2}\right)K(\v{x}),
\end{align*}
where $x_j$ is the $j$-th component of $\v{x}\in\mathbb{R}^q$.

Under appropriate conditions on the kernel $K$~\cite{stefanski1990deconvolving}, we can see that 
\begin{align} \label{eqn:unbiased_prob}
\mathbb{E}\{\widehat{\phi}_{\v{x}}(\v{x})|\{\v{x}_i\}_{i=1}^n\}=\widehat{\phi}^{\mathrm{np}}_{\v{x}}(\v{x}).
\end{align}
Therefore, $\widehat{\phi}_{\v{x}}(\v{x})$ in~\eqref{eqn:density_private} is effectively an unbiased estimate of $\widehat{\phi}^{\mathrm{np}}_{\v{x}}(\v{x})$ in~\eqref{eqn:density_no_privacy}. In average, we are canceling the effect of the differential privacy noise. \revise{Selecting bandwidth (or smoothing parameter) $h$ is an important aspect of kernel estimation. In~\cite{Parzen1962}, it was shown that $\lim_{n\rightarrow \infty} h=0$ guarantees asymptotic ubiasedness (i.e., point-wise convergence of the kernel density estimate to the true density function) while $\lim_{n\rightarrow \infty} nh=+\infty$ is required to ensure asymptotic consistency. Many studies have focused on finding optimal bandwidth~\cite{rudemo1982empirical,robert1976choice,loftsgaarden1965,koontz1972asymptotic}. Numerical methods based on cross validation for setting the bandwidth are proposed in~\cite{delaigle2002estimation,delaigle2004practical}. Often, it is recommended to compare the results from different bandwidth selection algorithms to avoid misleading conclusions caused by over-smoothing or under-smoothing of the density estimate~\cite{wang2019simple}.  These results have been also extended to noisy measurements with deconvoluting kernel density estimation~\cite{stefanski1990deconvolving,es2005asymptotic}. If} $h$ scales according to $n^{-1/5}$, $\widehat{\phi}_{\v{x}}(\v{x})$ is a consistent estimator of $\phi_{\v{x}}$ as $n\rightarrow\infty$, i.e., $\widehat{\phi}_{\v{x}}(\v{x})$ converges $\phi_{\v{x}}$ point-wise for all $\v{x}\in\mathrm{supp}(\phi_{\v{x}})$~\cite{stefanski1990deconvolving}. \revise{
Note that by selecting $h=\mathcal{O}(n^{-1/5})$, we get
\begin{align*}
\int \mathbb{E}\{\widehat{\phi}_{\v{x}}(\v{x})-\phi_{\v{x}}(\v{x})\}^2=\mathcal{O}(n^{-4/5}),
\end{align*}
where $\mathcal{O}$ denotes the Bachmann–Landau notation. Therefore, $\int \mathbb{E}\{\widehat{\phi}_{\v{x}}(\v{x})-\phi_{\v{x}}(\v{x})\}^2\rightarrow 0$  as $n\rightarrow \infty$. This means that the effect of the differential-privacy noise is effectively negligible on large datasets. 
}

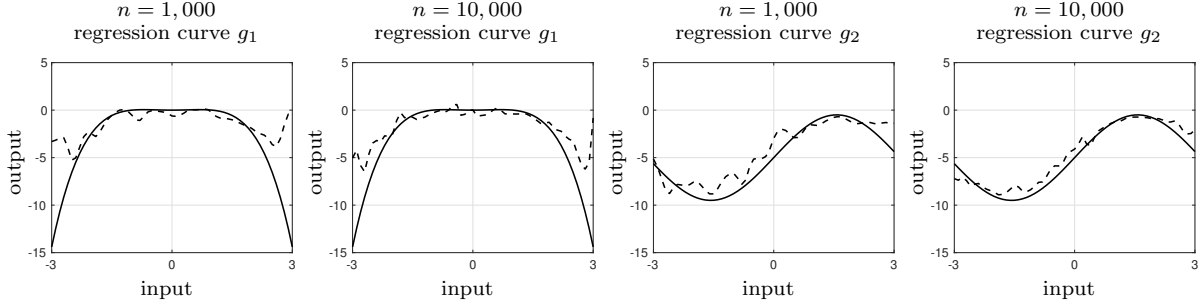
\begin{figure}
	\centering
	\begin{tikzpicture}
	\node[] at (0,0) {\includegraphics[width=.25\linewidth]{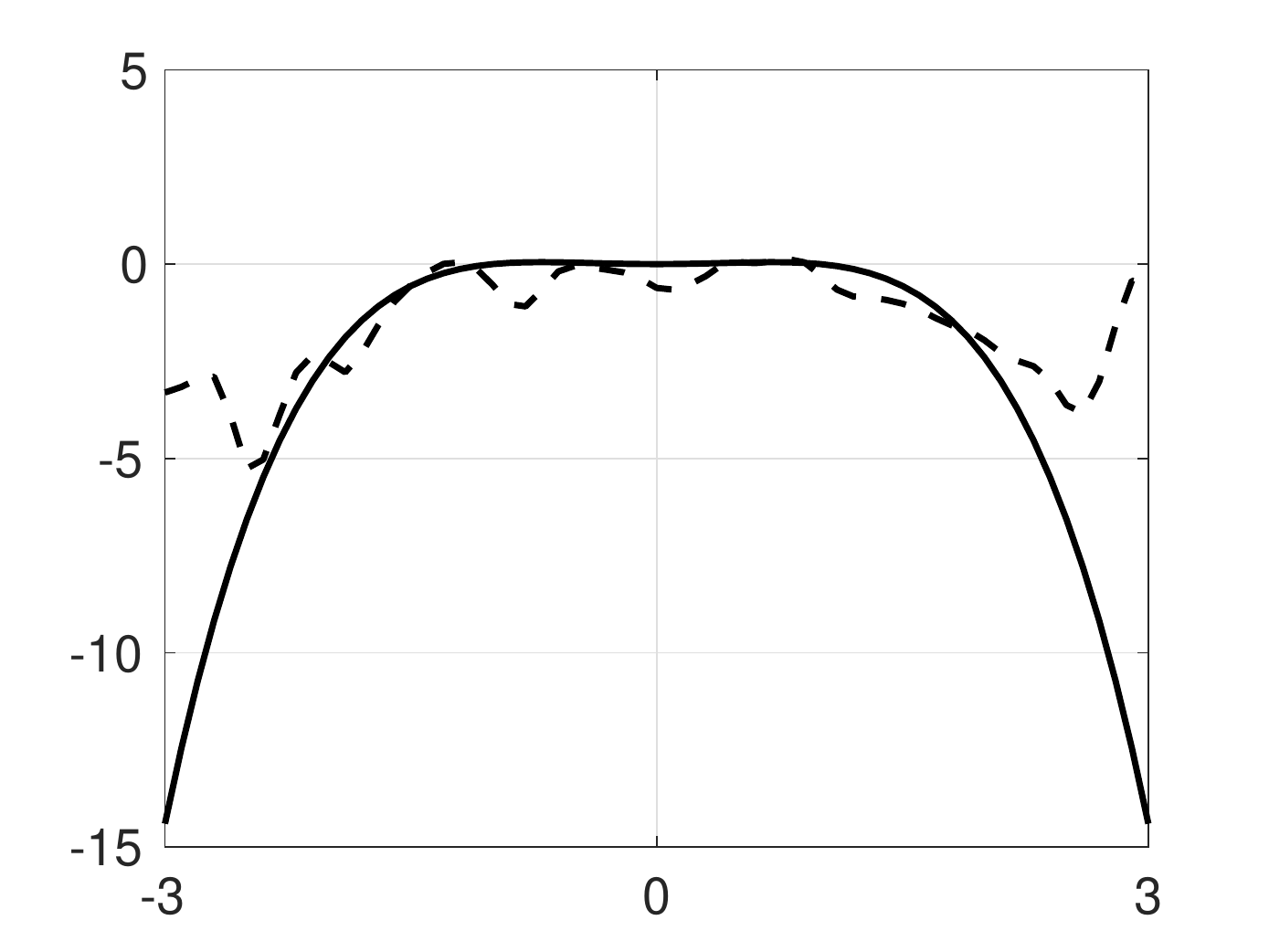}};
	\node[] at (4,0) {\includegraphics[width=.25\linewidth]{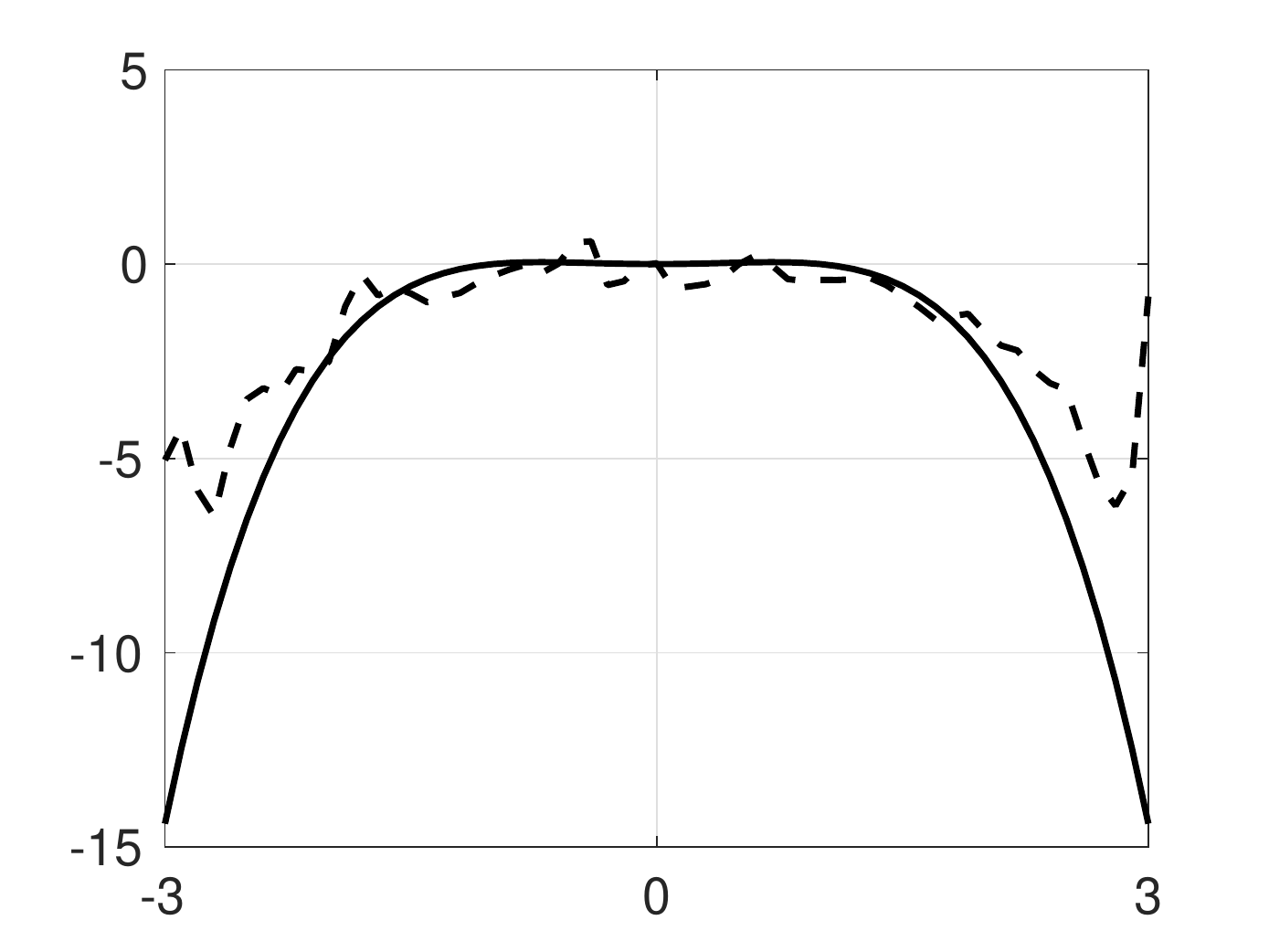}};
	\node[] at (8,0) {\includegraphics[width=.25\linewidth]{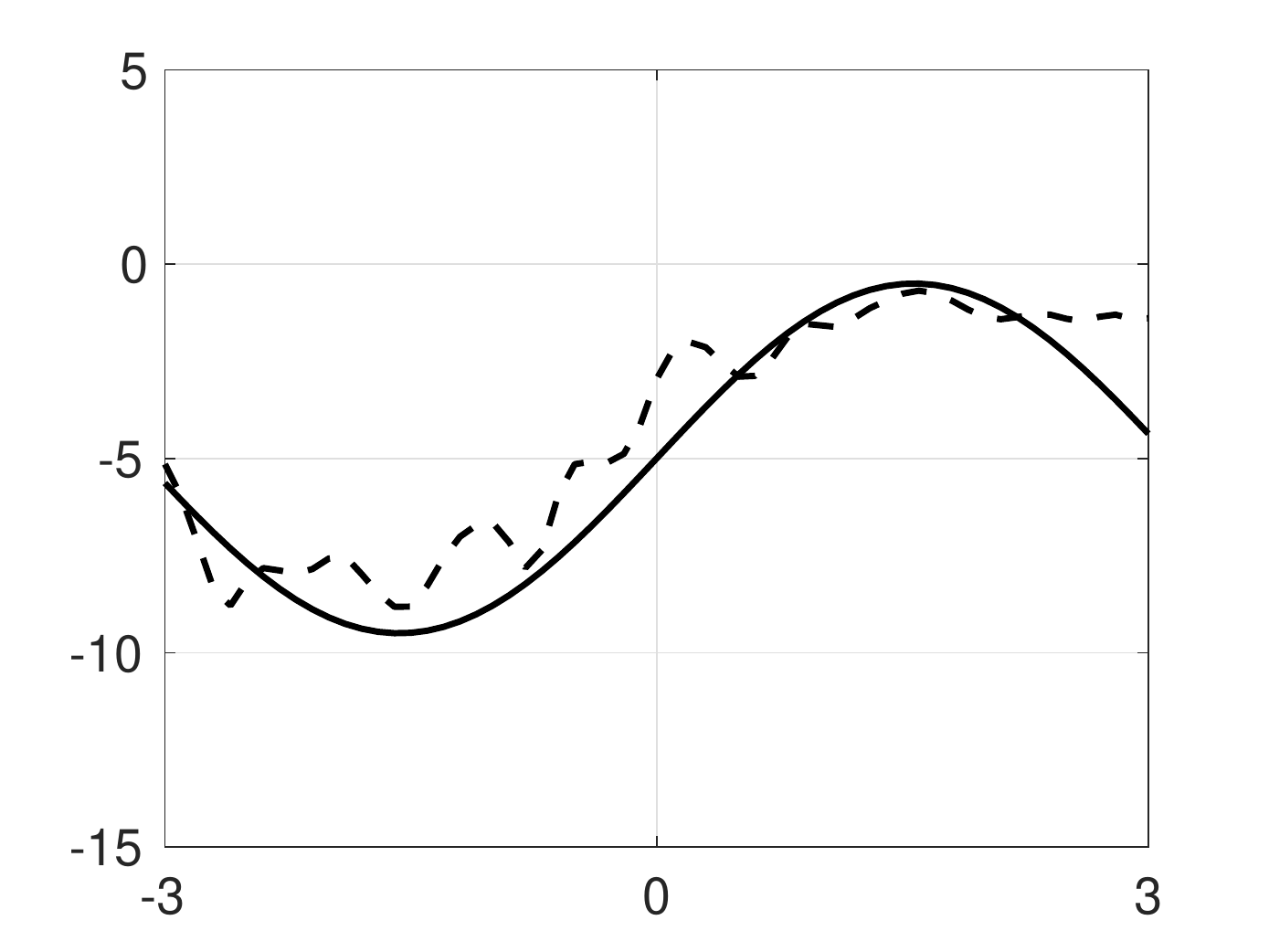}};
	\node[] at (12,0) {\includegraphics[width=.25\linewidth]{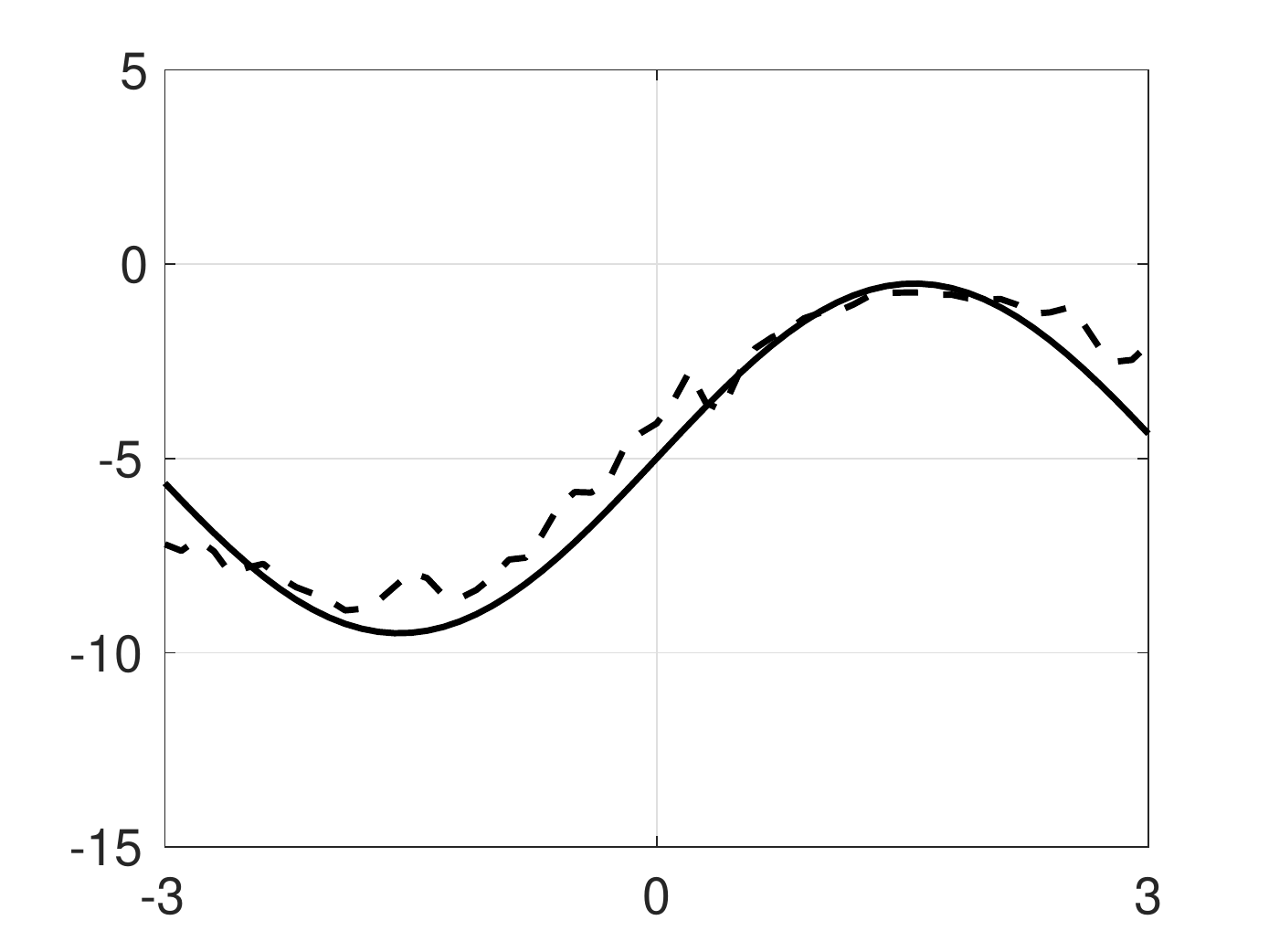}};
	\node[] at (0,-1.7) {\footnotesize input};
	\node[] at (4,-1.7) {\footnotesize input};
	\node[] at (8,-1.7) {\footnotesize input};
	\node[] at (12,-1.7) {\footnotesize input};
	\node[] at (0,+2) {\footnotesize $n=1,000$};
	\node[] at (0,+1.7) {\footnotesize regression curve $g_1$};
	\node[] at (4,+2) {\footnotesize $n=10,000$};
	\node[] at (4,+1.7) {\footnotesize regression curve $g_1$};
	\node[] at (8,+2) {\footnotesize $n=1,000$};
	\node[] at (8,+1.7) {\footnotesize regression curve $g_2$};
	\node[] at (12,+2) {\footnotesize $n=10,000$};
	\node[] at (12,+1.7) {\footnotesize regression curve $g_2$};
	\node[rotate=90] at (-2.,0) {\footnotesize output};
	\node[rotate=90] at (2.,0) {\footnotesize output};
	\node[rotate=90] at (6.,0) {\footnotesize output};
	\node[rotate=90] at (10.,0) {\footnotesize output};
	\end{tikzpicture}
	\caption{
		\label{fig:simulation_Gaussian_mix}
		The kernel regression model (dashed black) and true regression curve (solid black) for mixture Gaussian data made differentially private with $\epsilon=10$. 
	}
\end{figure}

\begin{figure}
	\centering
	\begin{tikzpicture}
	\node[] at (0,0) {\includegraphics[width=.25\linewidth]{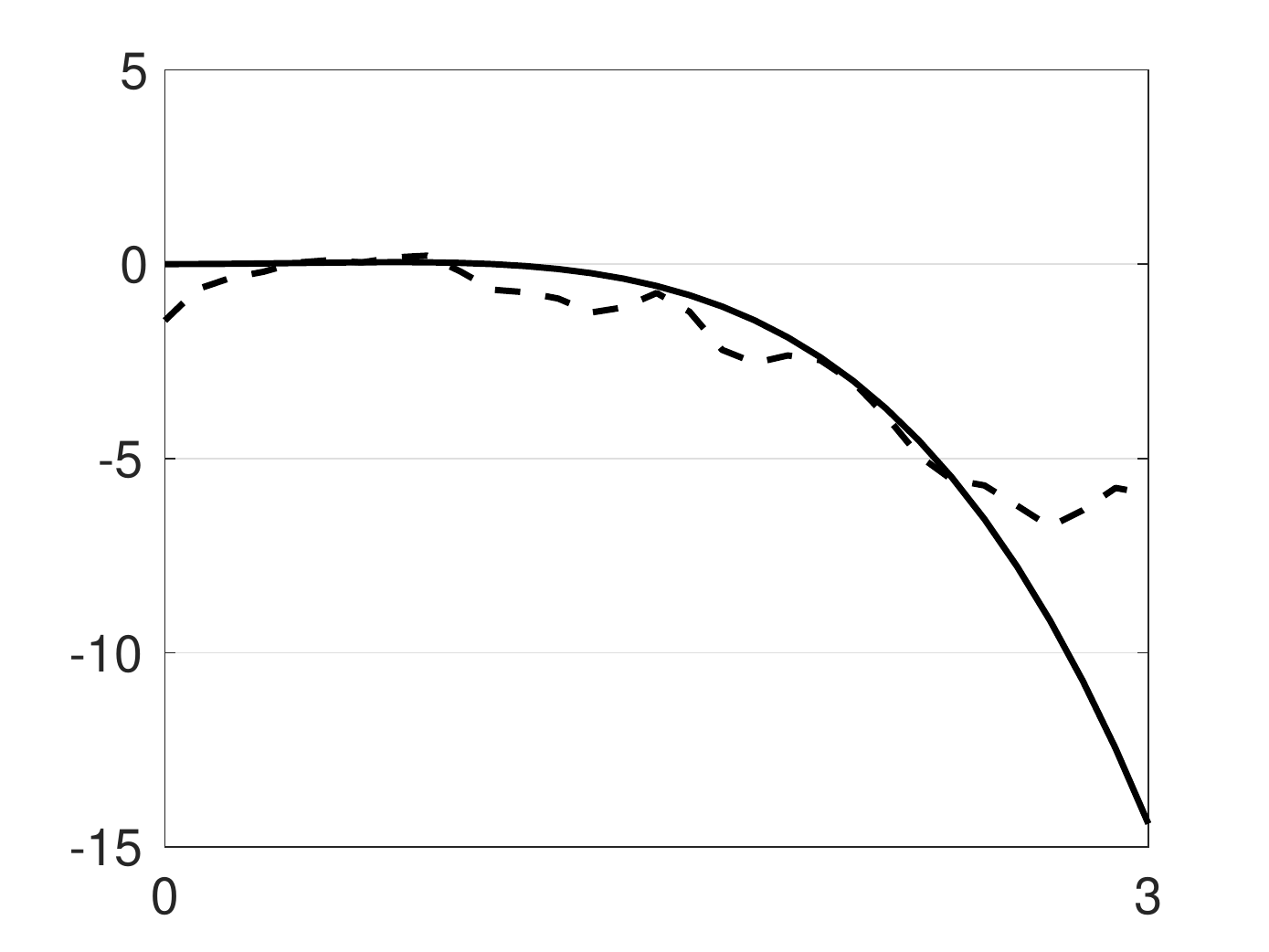}};
	\node[] at (4,0) {\includegraphics[width=.25\linewidth]{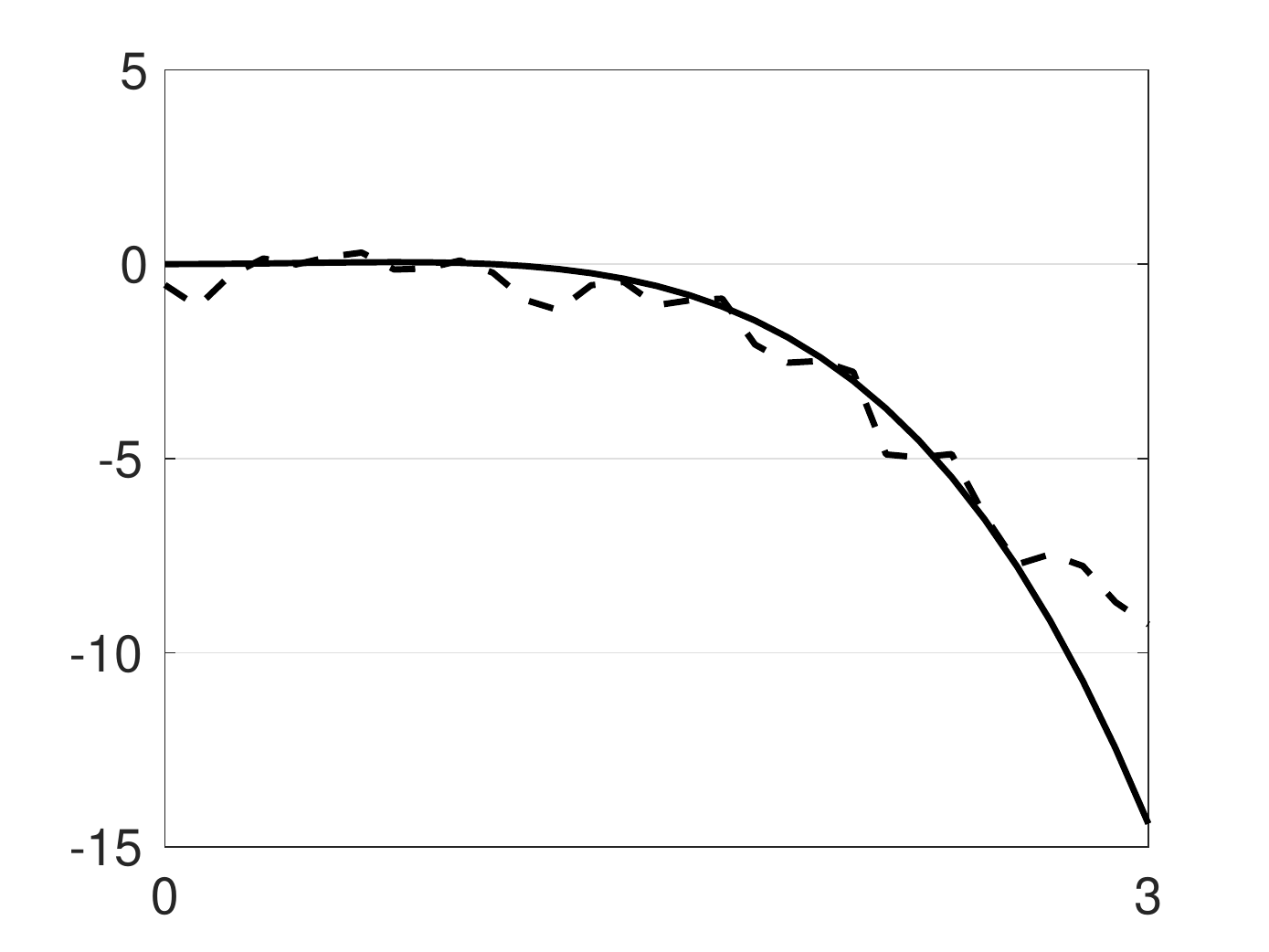}};
	\node[] at (8,0) {\includegraphics[width=.25\linewidth]{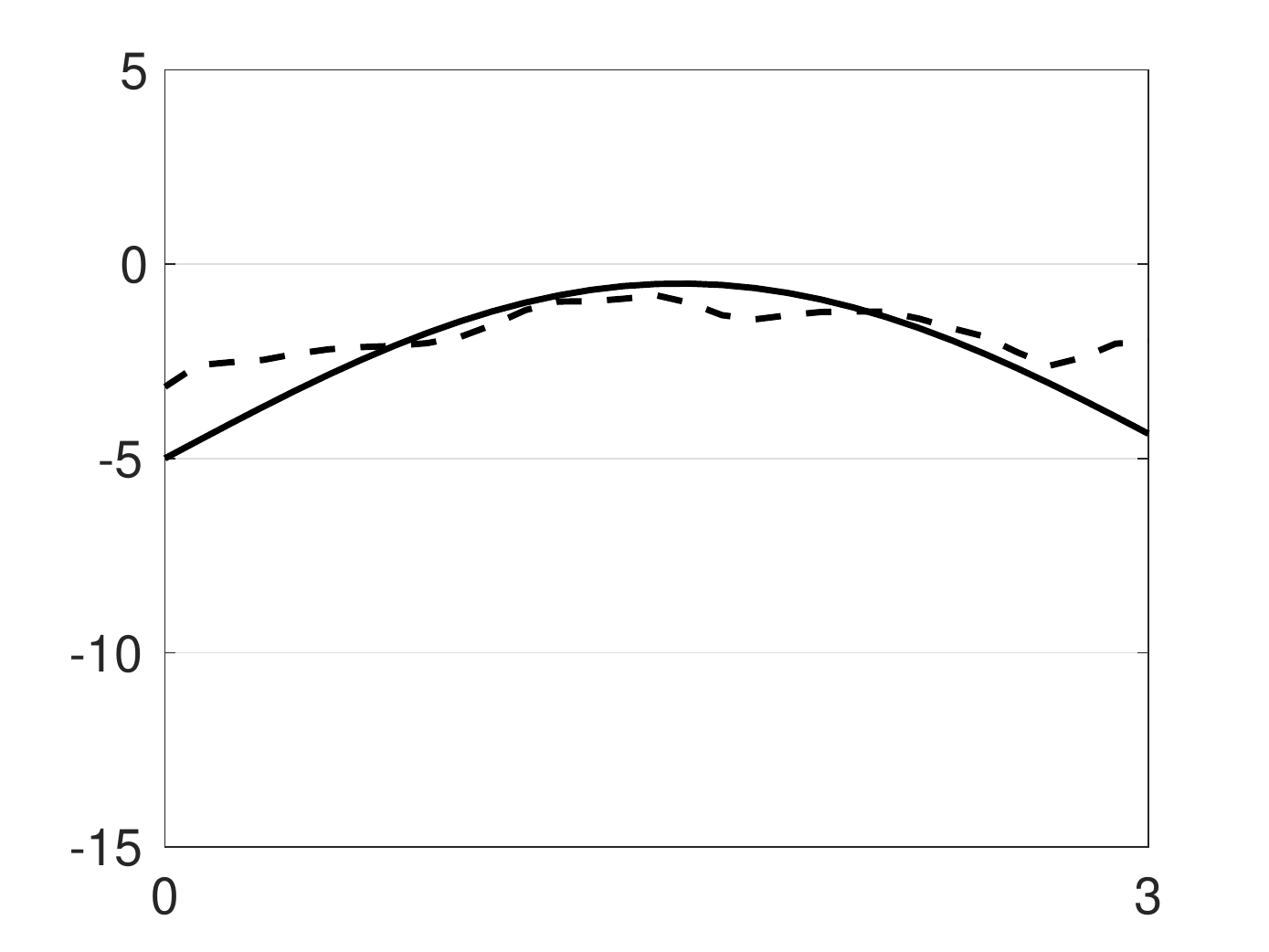}};
	\node[] at (12,0) {\includegraphics[width=.25\linewidth]{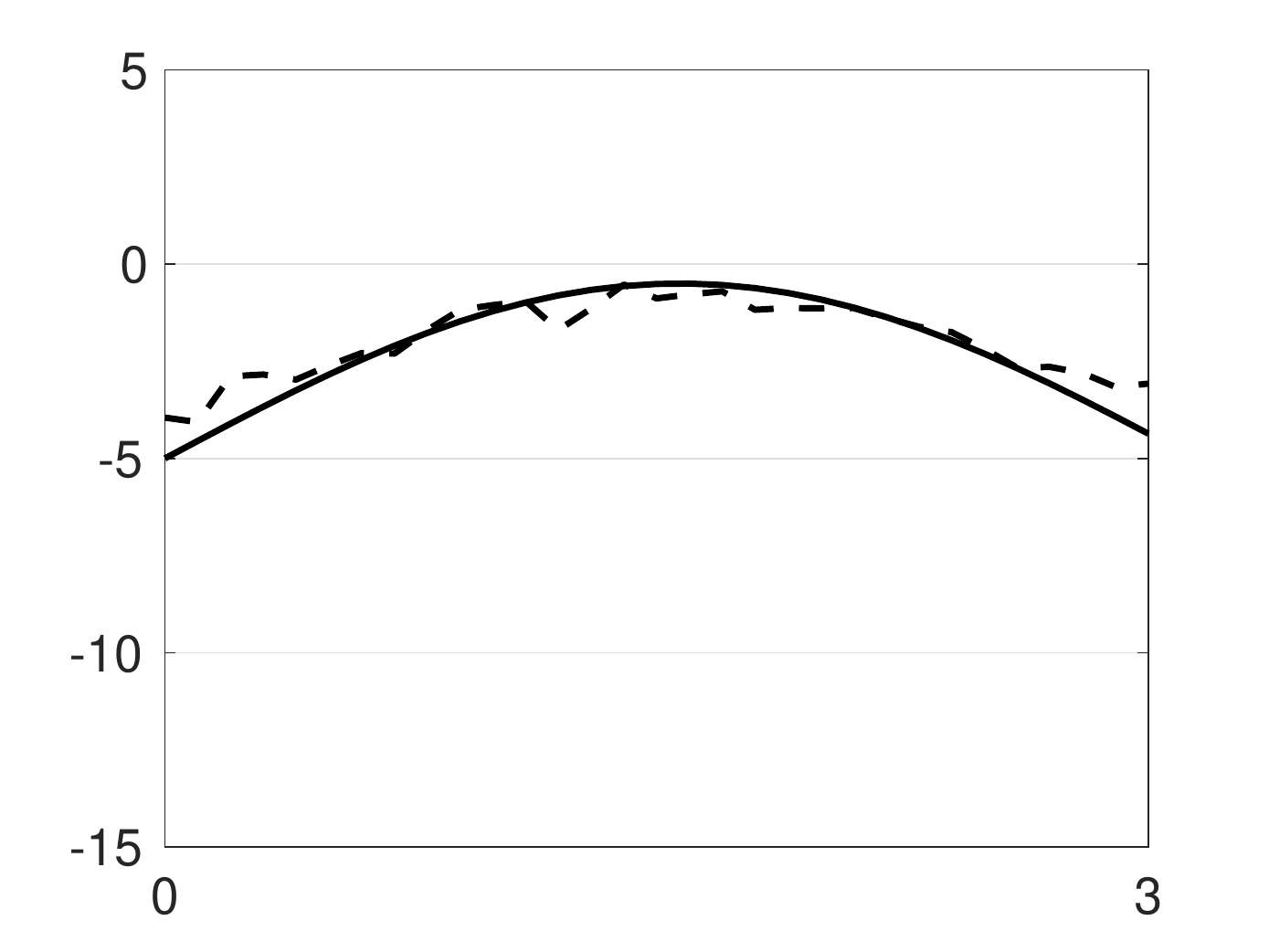}};
	\node[] at (0,-1.7) {\footnotesize input};
	\node[] at (4,-1.7) {\footnotesize input};
	\node[] at (8,-1.7) {\footnotesize input};
	\node[] at (12,-1.7) {\footnotesize input};
	\node[] at (0,+2) {\footnotesize $n=1,000$};
	\node[] at (0,+1.7) {\footnotesize regression curve $g_1$};
	\node[] at (4,+2) {\footnotesize $n=10,000$};
	\node[] at (4,+1.7) {\footnotesize regression curve $g_1$};
	\node[] at (8,+2) {\footnotesize $n=1,000$};
	\node[] at (8,+1.7) {\footnotesize regression curve $g_2$};
	\node[] at (12,+2) {\footnotesize $n=10,000$};
	\node[] at (12,+1.7) {\footnotesize regression curve $g_2$};
	\node[rotate=90] at (-2.,0) {\footnotesize output};
	\node[rotate=90] at (2.,0) {\footnotesize output};
	\node[rotate=90] at (6.,0) {\footnotesize output};
	\node[rotate=90] at (10.,0) {\footnotesize output};
	\end{tikzpicture}
	\caption{
		\label{fig:simulation_chi_squared}
		The kernel regression model (dashed black) and true regression curve (solid black) for chi-squared data made differentially private with $\epsilon=10$. 
	}
\end{figure}
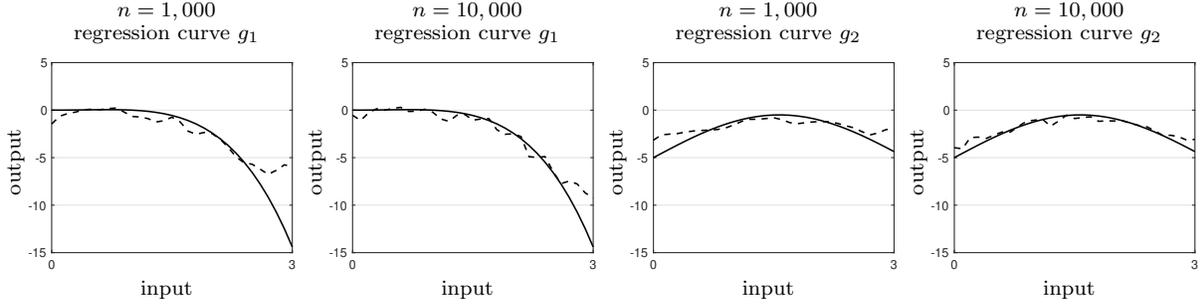

For regression analysis, we consider independently distributed data points $\{(\v{x}[i],\v{y}[i])\}_{i=1}^n$ from common probability density function. We would like to understand the relationship between inputs $\v{x}[i]$ and outputs $\v{y}[i]$ for all $i$. Similarly, we assume that we can only access noisy privacy-preserving inputs $\{\v{z}[i]\}_{i=1}^n$ instead of accurate inputs $\{\v{x}[i]\}_{i=1}^n$. Following the argument above, we can also construct the Nadaraya-Watson kernel regression (see, e.g.,~\cite{hardle1990applied}) as
\begin{align} \label{eqn:kernel_regression}
\widehat{m}(\v{x}):=\frac{\sum_{i=1}^n \widehat{K}_h((\v{x}-\v{z}[i])/h)\v{y}[i]}{\sum_{i=1}^n \widehat{K}_h((\v{x}-\v{z}[i])/h)}.
\end{align}
Under appropriate conditions on the kernel $K$ and the bandwidth $h$~\cite{fan1993nonparametric}, $\widehat{m}(\v{x})$ converges to $\mathbb{E}\{\v{y}|\v{x}\}$ almost surely. In practice the bandwidth can be computed by minimizing the cross-validation cost, i.e., the error of estimating each $\v{y}[j]$ using  the Nadaraya-Watson kernel regression constructed from  $\{(\v{z}[i],\v{y}[i])\}_{i\in\{1,\dots,n\}\setminus\{j\}}$ averaged over all choices of $\ell$. \revise{The optimal bandwidth is given by
\begin{align}\label{eqn:CV}
\argmin_{h} \sum_{j=1}^n \ell(\v{y}[j],\widehat{m}_{-j}(\v{x}[j])),
\end{align}
where $\ell$ is a fitness function, e.g., $\ell(\v{y},\v{y}')=\|\v{y}-\v{y}'\|_2^2$, and $\widehat{m}_{-j}(\v{x})$ is the Nadaraya-Watson kernel regression constructed from $\{(\v{z}[i],\v{y}[i])\}_{i\in\{1,\dots,n\}\setminus\{j\}}$:
\begin{align*}
\widehat{m}_{-j}(\v{x}):=\frac{\sum_{i\in\{1,\dots,n\}\setminus\{j\}} \widehat{K}_h((\v{x}-\v{z}[i])/h)\v{y}[i]}{\sum_{i\in\{1,\dots,n\}\setminus\{j\}} \widehat{K}_h((\v{x}-\v{z}[i])/h)}.
\end{align*}
This approach has been widely used for setting the bandwidth in non-parametric regression~\cite{delaigle2007nonparametric}.}

\section*{Results}
In this section, we demonstrate the performance of the developed methods on \revise{multiple datasets. We first use a synthetic dataset for illustration purposes and then utilize real} financial and demographic datasets. \revise{Throughout this section, we use the following original kernel:
	\begin{align*}
	K(x)=\frac{1}{\pi} \frac{1}{1+x^2}.
	\end{align*}
	Note that $\v{x}=x$ is a scalar as we are only considering credit score as an input. This is the Cauchy distribution. We get the adjusted kernel in
	\begin{align*}
	\widehat{K}_h(x)
	&=\left(1-\frac{b^2}{h^2} \frac{\mathrm{d}^2}{\mathrm{d} x}\right)K(x)\\
	&=\frac{1}{\pi}\left[ \frac{1}{1+x^2} - \frac{b^2}{h^2}\frac{8x^2}{(x^2 + 1)^3} +\frac{b^2}{h^2} \frac{2}{(x^2 + 1)^2}\right].
	\end{align*}
	We use the cross-validation procedure in~\eqref{eqn:CV} to find the bandwidth in the following simulation and experiments. }

\begin{figure}[t]
	\centering
	\begin{tikzpicture}
	\node[] at (0,0) {\includegraphics[width=.45\linewidth]{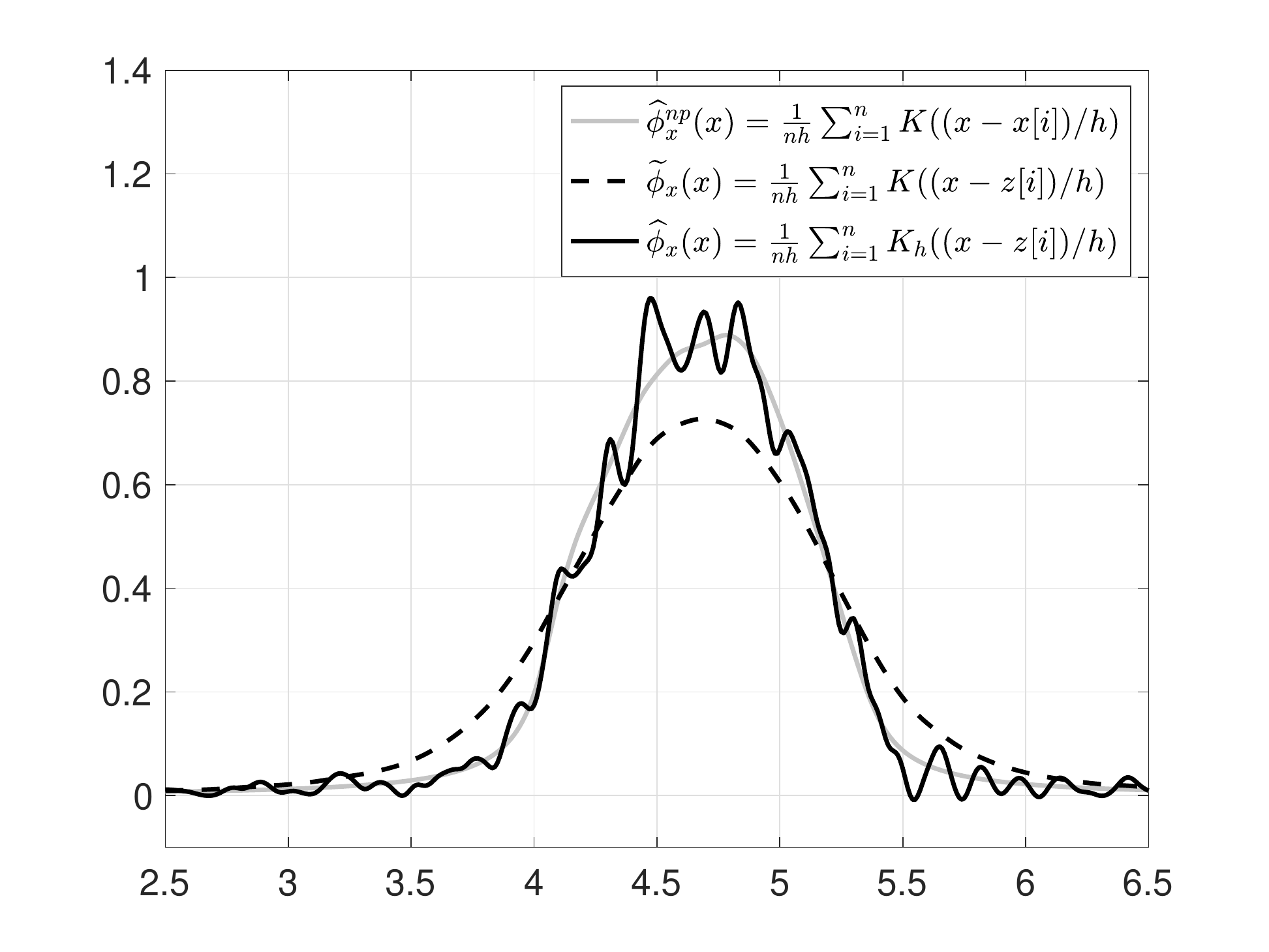}};
	\node[] at (0,-2.7) {\footnotesize $\log(\mbox{credit rating}-600)$};
	\node[rotate=90] at (-3.3,0) {\footnotesize probability density function};
	\end{tikzpicture}
	\caption{
		\label{fig:prob_fico}
		Estimates of probability density function of the credit score using original noiseless data with original kernel $\widehat{\phi}^{np}_{x}(x)=\frac{1}{nh}\sum_{i=1}^{n}K((x-x[i])/h)$ (solid gray), $\epsilon$-locally differential private data with original kernel
		$\widetilde{\phi}_{x}(x)=\frac{1}{nh}\sum_{i=1}^{n}K((x-z[i])/h)$ (dashed black), and $\epsilon$-locally differential private data with adjusted kernel
		$\widehat{\phi}_{x}(x)=\frac{1}{nh}\sum_{i=1}^{n}K_h((x-z[i])/h)$ (solid black) for $\epsilon=5.0$ and bandwidth $h=0.1$.
	}
\end{figure}

\begin{figure}[t]
	\centering
\begin{tikzpicture}
\node[] at (0,0) {\includegraphics[width=.45\linewidth]{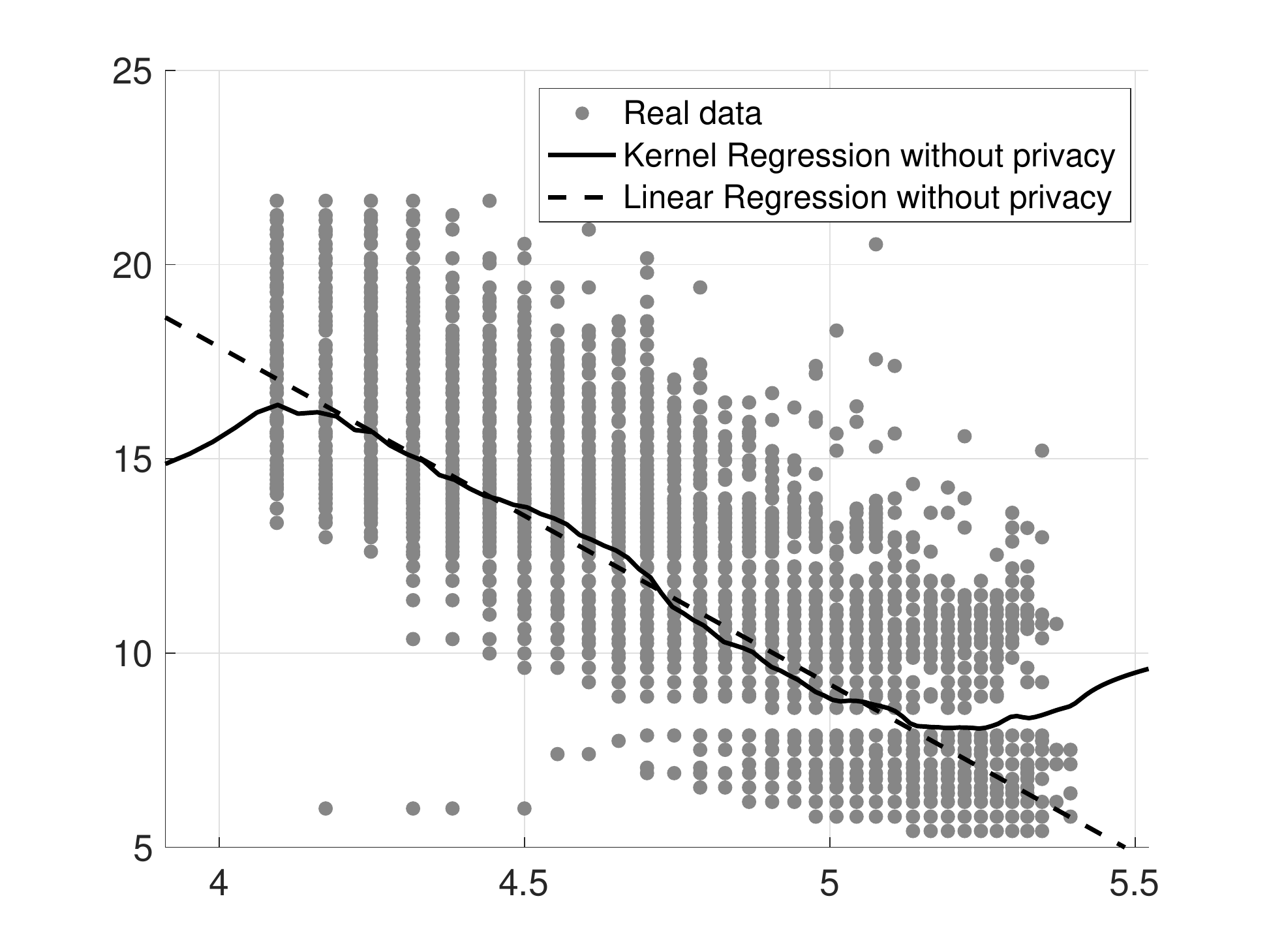}};
\node[] at (0,-2.7) {\footnotesize $\log(\mbox{credit rating}-600)$};
\node[rotate=90] at (-3.3,0) {\footnotesize interest rate (percentage)};
\end{tikzpicture}
\caption{
\label{fig:no_noise_linear_kernel}
	The kernel regression model (solid black) and the linear regression model (dashed black)  based on the original data with bandwidth $h=0.02$ superimposed on the original noiseless data (gray dots). The mean squared error for the kernel regression model is $4.42$ and the mean squared error for the linear regression model is $4.61$. 
}
	\centering
\begin{tikzpicture}
\node[] at (0,0) {\includegraphics[width=.45\linewidth]{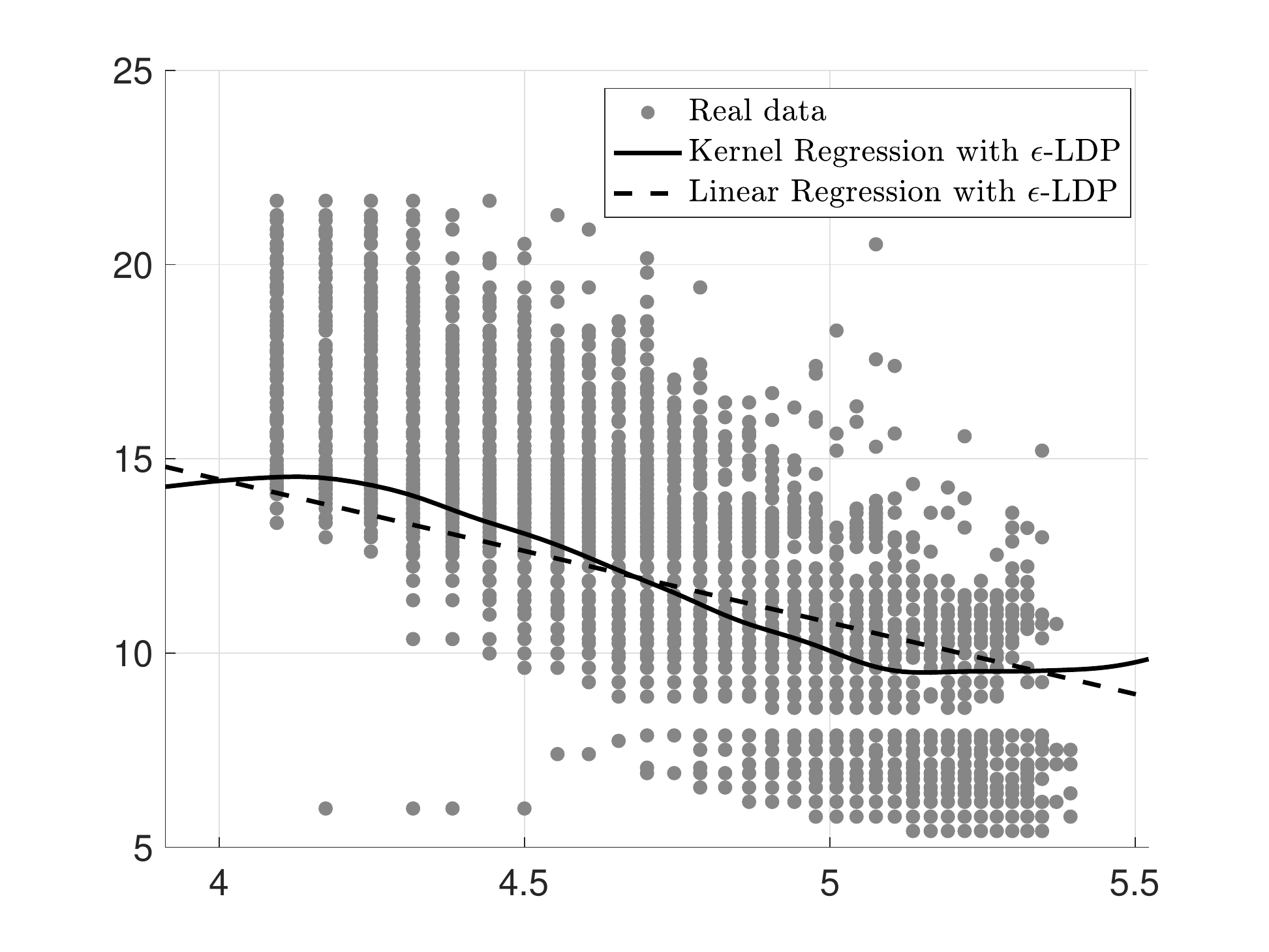}};
\node[] at (0,-2.7) {\footnotesize $\log(\mbox{credit rating}-600)$};
\node[rotate=90] at (-3.3,0) {\footnotesize interest rate (percentage)};
\end{tikzpicture}
	\caption{
		\label{fig:noisy_epsilon_5_linear_kernel}
		The kernel regression model (solid black) and the linear regression model (dashed black)  based on the $\epsilon$-locally differential private data with $\epsilon=5$ and bandwidth $h=0.20$ superimposed on the original noiseless data (gray dots). The mean squared error for the kernel regression model is $5.70$ and the mean squared error for the linear regression model is $7.11$. 
	}
\end{figure}

\begin{figure}[t]
	\centering
\begin{tikzpicture}
\node[] at (0,0) {\includegraphics[width=.45\linewidth]{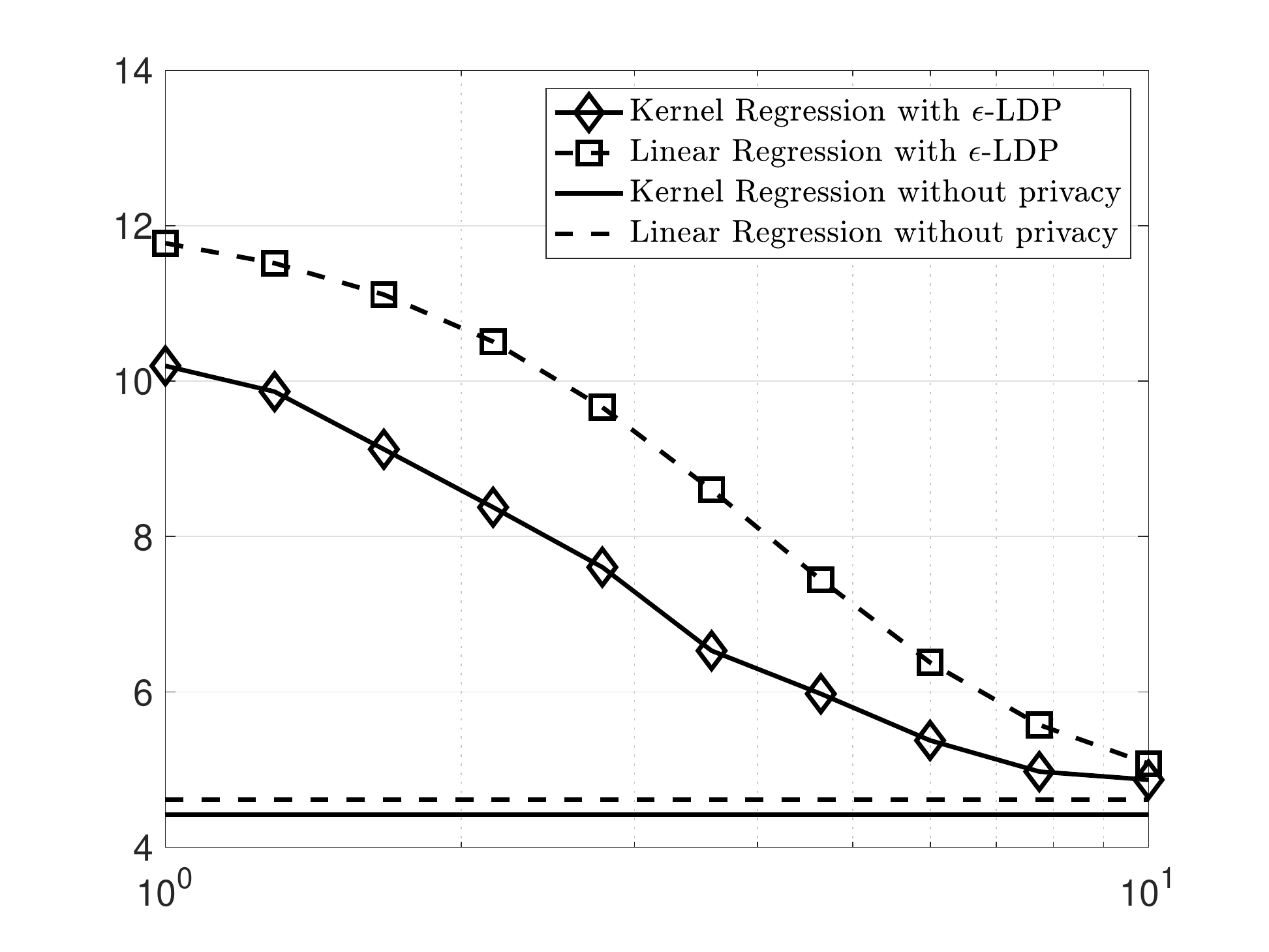}};
\node[] at (0,-2.7) {\footnotesize privacy budget};
\node[rotate=90] at (-3.2,0) {\footnotesize mean squared error};
\end{tikzpicture}
	\caption{
		\label{fig:MSE_vs_epsilon}
		The mean squared error for the kernel regression model and the linear regression model based on the $\epsilon$-locally differential private ($\epsilon$-LDP in the legend) data versus privacy budget $\epsilon$. The horizontal lines show the mean squared error for the kernel regression model and the linear regression model based on original noiseless data. 
	}
\end{figure}

\revise{
	\subsection*{Synthetic Dataset}
	We use a simulation study is to illustrate the performance of the Nadaraya-Watson kernel regression in~\eqref{eqn:kernel_regression} for privacy-preserving data. We consider multiple scenarios. We use two distributions for $\{x[i]\}_{i=1}^n$.  The first one is a Gaussian mixture $(1/3)\mathcal{N}(-1,1)+(2/3)\mathcal{N}(3/2,1/2)$ truncated over $[-3,3]$. The second distribution is a chi-squared distribution
	with three degrees of freedom distribution $\chi^2(3)$ truncated over $[0,3]$. The truncation in both cases is for satisfaction of Assumption~\ref{assum:1}.  We also consider two regression curves: $g_1:x\mapsto x^2(1-x^2)/5$ and $g_2:x\mapsto 4.5\sin(x)-5$. Finally, we assume a Gaussian measurement noise $\mathcal{N}(0,1)$, i.e., $y[i]=g_j(x[i])+v[i]$ for $j=1,2$, where $v[i]$ is a zero mean Gaussian random variable with unit variance. 
	
	Figure~\ref{fig:simulation_Gaussian_mix} shows the kernel regression model (dashed black) and true regression curve (solid black) for mixture Gaussian data made differentially private with $\epsilon=10$. Here, we consider two dataset size of $n=1,000$ and $n=10,000$ and two regression curves of $g_1$ and $g_2$, introduced earlier.  The Nadaraya-Watson kernel regression using differentially-private provides fairly accurate predictions. The accuracy of the prediction improves as the dataset gets larger. Figure~\ref{fig:simulation_chi_squared} illustrates the kernel regression model (dashed black) and true regression curve (solid black) for chi-squared data made differentially private with $\epsilon=10$. This shows that the fitness of the Nadaraya-Watson kernel regression is somewhat independent of the underlying distribution of the data.
	
}

\subsection*{Lending Club Dataset}
The dataset contains information of 2,260,701 accepted and 27,648,741 rejected loans application on Lending Club, a peer-to-peer lending platform, over 2007 to 2018. The dataset is available for download on Kaggle~\cite{kaggle1}. For the accepted loans, dataset contains interest rates of the loans per annum and loan attributes, such as total loan size, and borrower information, such as number of credit lines, credit rating, state of residence, and age. Here, we only focus on data from 2010 (to avoid possible yearly fluctuations of the interest rate), which contains 12,537 accepted loans. We also focus on the relationship between the FICO\footnote{\url{https://www.fico.com/en/products/fico-score}} credit score (low range) and the interest rates of the loan. This is an interesting relationship  
pointing to the value of credit rating reports~\cite{czarnitzki2007credit}. The FICO credit score is very sensitive (as it relates to the financial health of an individual) and possesses a significant commercial value (as it is sold by a for-profit corporation). Thus, we assume that is is made available publicly in a privacy-preserving manner using~\eqref{eqn:additive}. Note that the original data in~\cite{kaggle1} provides this data in an anonymized manner without privacy-preserving noise.

Figure~\ref{fig:prob_fico} illustrates estimates of probability density function of the credit score $\phi_x(x)$ using original noiseless data with original kernel $\widehat{\phi}^{np}_{x}(x)$ in~\eqref{eqn:density_no_privacy} (solid gray), $\epsilon$-locally differential private data with original kernel
$\widetilde{\phi}_{x}(x)=\frac{1}{nh}\sum_{i=1}^{n}K((x-z[i])/h)$ (dashed black), and $\epsilon$-locally differential private data with adjusted kernel
in~\eqref{eqn:density_private} (solid black) for $\epsilon=5.0$ and bandwidth $h=0.1$. Note that $\widetilde{\phi}_{x}(x)=\frac{1}{nh}\sum_{i=1}^{n}K((x-z[i])/h)$ is a naive density estimate as it does not try to cancel the effect of the privacy-preserving noise. Clearly, using the original kernel for the noisy privacy-preserving data flattens the density estimate $\widetilde{\phi}_{x}(x)$. This is because we are in fact observing a convolution of the original probability density with the probability density of the Laplace noise. Upon using the adjusted kernel $\widehat{K}_h(x)$ the estimate of the probability density using the noisy privacy-preserving data matches the estimate of the probability density with the original data (with additional fluctuations due to the presence of noise). This provides a numerical validation of~\eqref{eqn:unbiased_prob}.

Now, let us focus on the regression analysis. Figure~\ref{fig:no_noise_linear_kernel} shows the kernel regression model (solid black) and the linear regression model (dashed black)  based on the original data with bandwidth $h=0.02$ superimposed on the original noiseless data (gray dots). The mean squared error for the kernel regression model is $4.42$ and the mean squared error for the linear regression model is $4.61$. The kernel regression model is thus slightly superior (roughly 4\%) to the linear regression model; however, the gap is narrow. Figure~\ref{fig:noisy_epsilon_5_linear_kernel} illustrates the kernel regression model (solid black) and the linear regression model (dashed black) based on the $\epsilon$-locally differential private data with $\epsilon=5$ and bandwidth $h=0.20$ superimposed on the original noiseless data (gray dots). The mean squared error for the kernel regression model is $5.70$ and the mean squared error for the linear regression model is $7.11$. In this case, the kernel regression model is considerably (roughly 20\%) better. In Figure~\ref{fig:MSE_vs_epsilon}, we observe the mean squared error for the kernel regression model and the linear regression model based on the $\epsilon$-locally differential private data versus privacy budget $\epsilon$. Clearly, the kernel regression model is consistently superior to the linear regression model. As $\epsilon$ grows larger, the performance of the kernel regression model and the linear regression model based on the $\epsilon$-locally differential private data converge to the performance of the kernel regression model and the linear regression model based on original noiseless data. This intuitively makes sense as, by increasing the privacy budget, the magnitude of the privacy-preserving noise becomes smaller.

\begin{figure}[t]
	\centering
\begin{tikzpicture}
\node[] at (0,0) {\includegraphics[width=.45\linewidth]{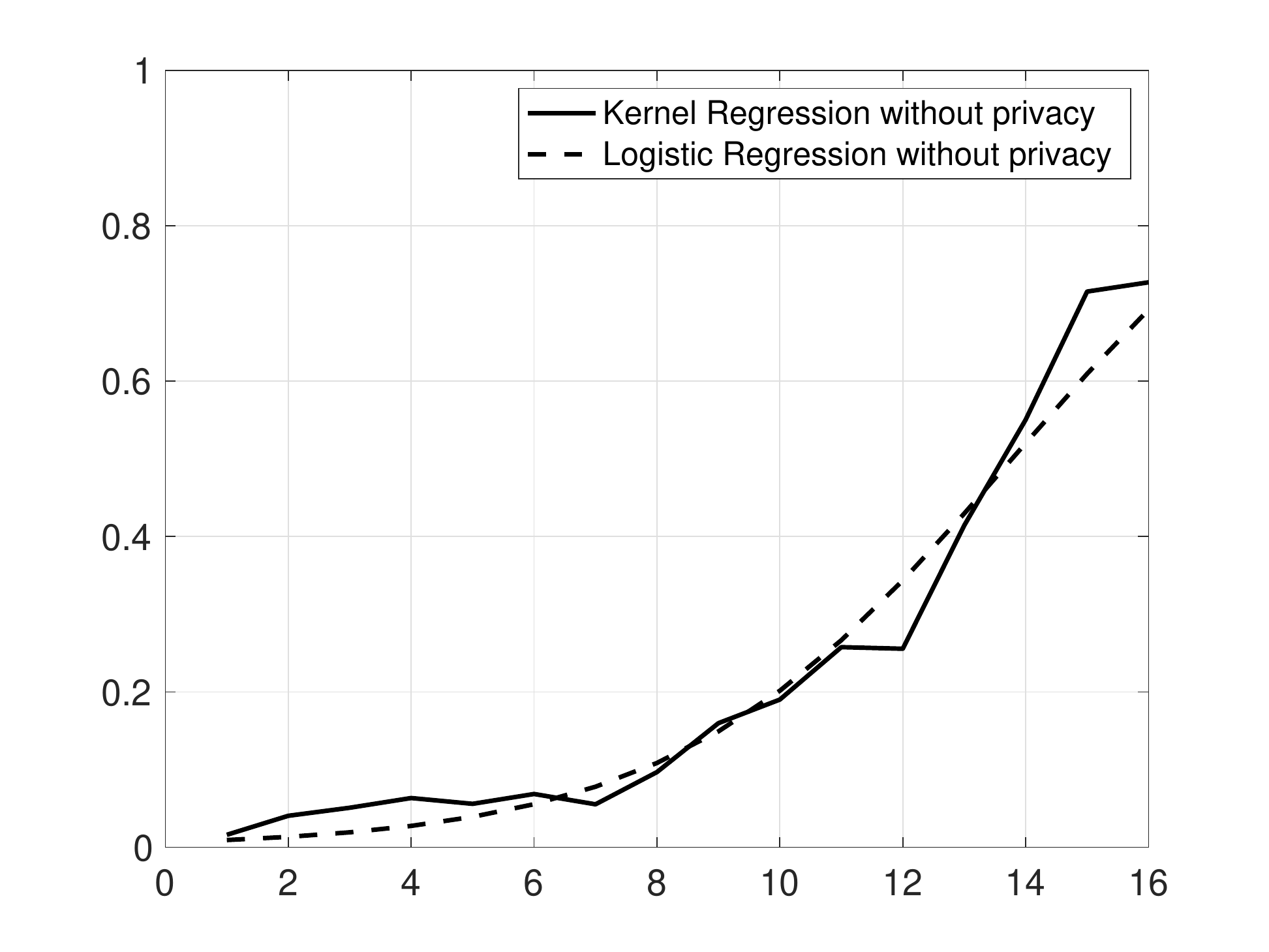}};
\node[] at (0,-2.7) {\footnotesize education (years)};
\node[rotate=90] at (-3.4,0) {\footnotesize $\mathbb{P}\{\mbox{income}\geq \mbox{50,000\$}\}$};
\end{tikzpicture}
	\caption{
		\label{fig:adult_no_noise_linear_kernel}	
		The kernel regression model (solid black) and the logistic regression model (dashed black)  based on the original data with bandwidth $h=0.17$. The logarithm of the likelihood for the kernel regression model is $-0.49$ and the logarithm of the likelihood for the logistic regression model is $-0.50$. 
	}
	\centering
\begin{tikzpicture}
\node[] at (0,0) {\includegraphics[width=.45\linewidth]{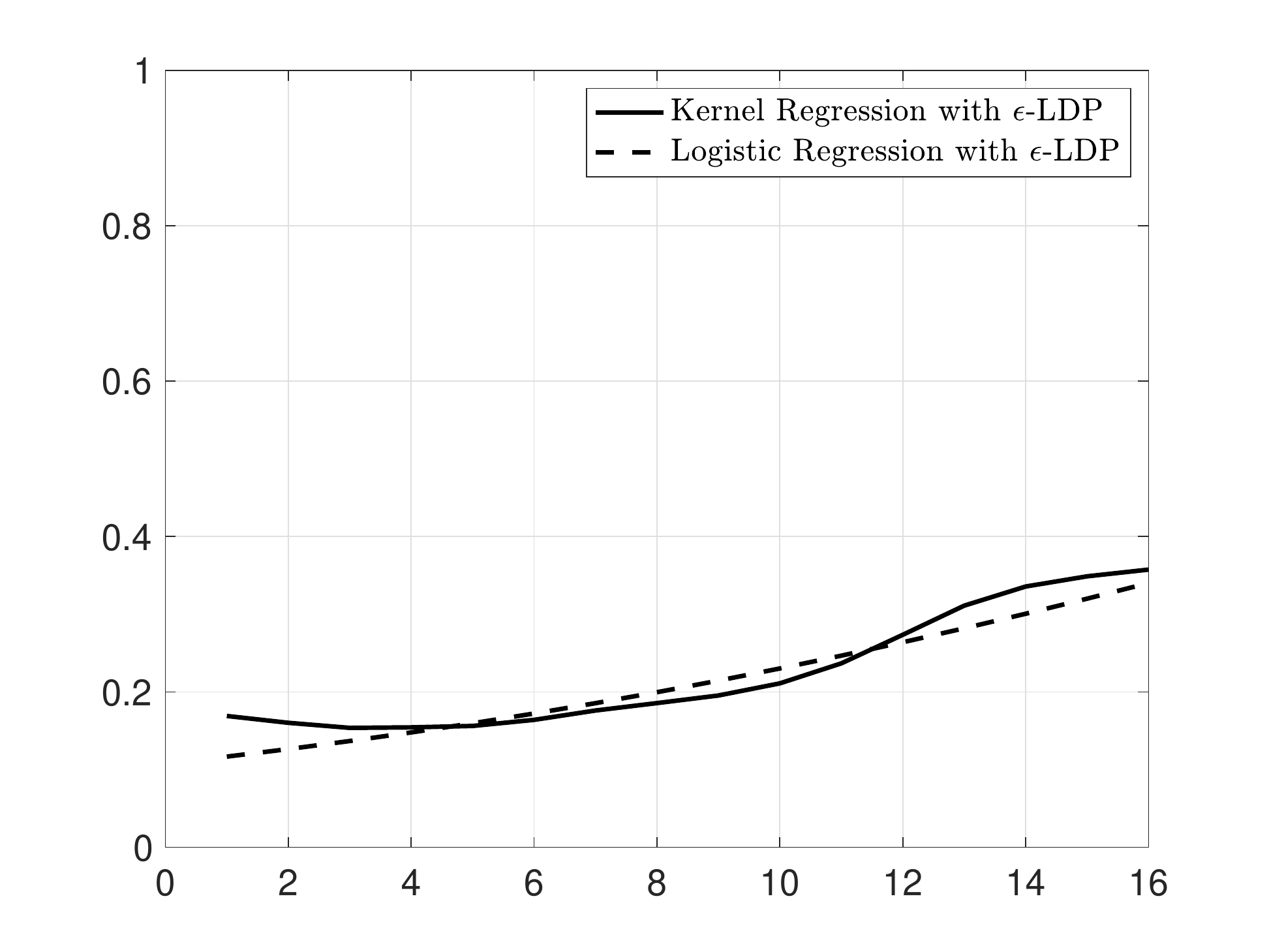}};
\node[] at (0,-2.7) {\footnotesize education (years)};
\node[rotate=90] at (-3.4,0) {\footnotesize $\mathbb{P}\{\mbox{income}\geq \mbox{50,000\$}\}$};
\end{tikzpicture}
	\caption{
		\label{fig:adult_noisy_epsilon_5_linear_kernel}
		The kernel regression model (solid black) and the logistic regression model (dashed black) based on the $\epsilon$-locally differential private data with $\epsilon=5.0$ bandwidth $h=2.98$. The logarithm of the likelihood for the kernel regression model is $-0.51$ and the logarithm of the likelihood for the logistic regression model is $-0.53$. 		
	}
\end{figure}

\begin{figure}[t]
	\centering
\begin{tikzpicture}
\node[] at (0,0) {\includegraphics[width=.45\linewidth]{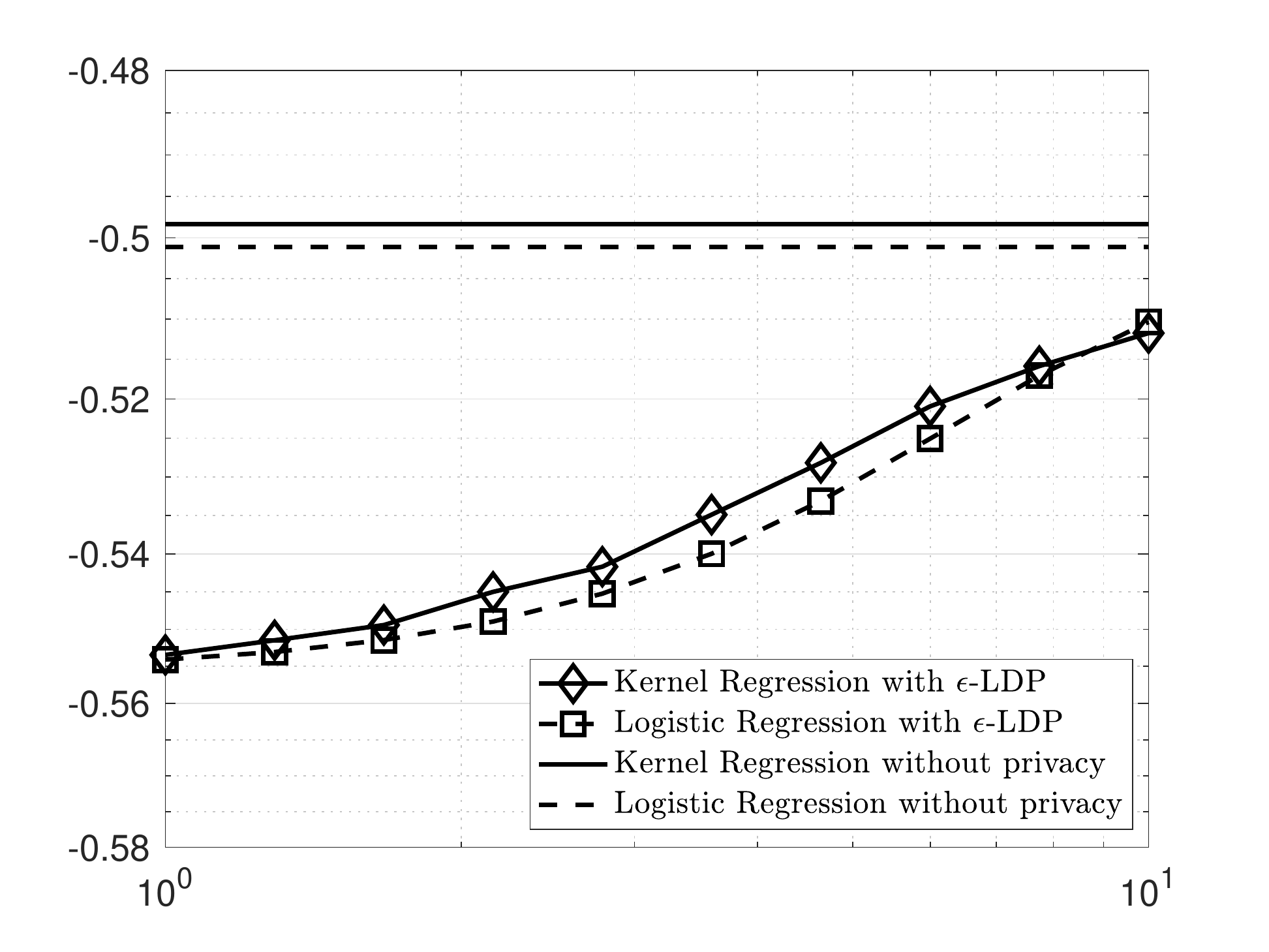}};
\node[] at (0,-2.7) {\footnotesize privacy budget};
\node[rotate=90] at (-3.5,0) {\footnotesize log likelihood};
\end{tikzpicture}
	\caption{
		\label{fig:adult_MSE_vs_epsilon}
		The logarithm of the likelihood for the kernel regression model and the logistic regression model based on the $\epsilon$-locally differential private data versus privacy budget $\epsilon$. The horizontal lines show the logarithm of the likelihood for the kernel regression model and the logistic regression model based on original noiseless data. 
	}
\end{figure}

\subsection*{Adult Dataset}
The dataset contains information of 32,561 individuals from the 1994 Census database. The dataset is available for download on UCI~\cite{Dua2019}. The dataset contains  attributes, such as education, age, work type, gender, race, and a binary report whether the individual earns more than 50,000\$ per year. We also focus on the relationship between the education (in years) and the individual ability to earn more than 50,000\$ per year. The education is  assumed to be made public in a privacy-preserving form following~\eqref{eqn:additive}. This information can be considered private as it can be used in conjunction with other information to de-anonymize the dataset.

Figure~\ref{fig:adult_no_noise_linear_kernel} The kernel regression model (solid black) and the logistic regression model (dashed black)  based on the original data with bandwidth $h=0.17$. The logarithm of the likelihood for the kernel regression model is $-0.49$ and the logarithm of the likelihood for the logistic regression model is $-0.50$. The kernel regression model is thus slightly superior (roughly 2\%) to the logistic regression model; however, the gap is almost negligible. Figure~\ref{fig:adult_noisy_epsilon_5_linear_kernel} illustrates the kernel regression model (solid black) and the logistic regression model (dashed black) based on the $\epsilon$-locally differential private data with $\epsilon=5.0$ bandwidth $h=2.98$. The logarithm of the likelihood for the kernel regression model is $-0.51$ and the logarithm of the likelihood for the logistic regression model is $-0.53$.  In this case, the kernel regression model is slightly (roughly 4\%) better. In Figure~\ref{fig:adult_MSE_vs_epsilon}, we observe the  logarithm of the likelihood for the kernel regression model and the logistic regression model based on the $\epsilon$-locally differential private data versus privacy budget $\epsilon$. The horizontal lines show the logarithm of the likelihood for the kernel regression model and the logistic regression model based on original noiseless data. Again, the kernel regression model is consistently superior to the logistic regression model. However, the effect is not as pronounced as the linear regression in the previous subsection. Finally, again, as $\epsilon$ grows larger, the performance of the kernel regression model and the logistic regression model based on the $\epsilon$-locally differential private data converge to the performance of the kernel regression model and the linear regression model based on original noiseless data. 

\section*{Discussion}

The density of privacy-preserving data is always flatter in comparison with the density function of the original data points due to convolution with privacy-preserving noise density function. This is certainly a cause for concern due to addition of differential-privacy noise in 2020 US Census. This unfortunate effect is always present irrespective of how many samples we gather  because we observe the convolution of the original probability density with the probability density of the privacy-preserving noise. This can result in miss-estimation of the heavy-hitters that often play an important role in social sciences due to their ties to minority groups. We developed density estimation methods  using smoothing kernels and used the framework of deconvoluting kernel density estimators  to remove the effect of privacy-preserving noise. This can result in a superior performance both for estimating probability density functions and for kernel regression in comparison to popular regression techniques, such as linear and logistic regression  models. In the case of estimating the probability density function, we could entirely remove the flatting effect of the privacy-preserving noise  at the cost of additional fluctuations. The fluctuations however could be reduced by gathering more data.

\bibliographystyle{ieeetr}
\bibliography{citation}

\section*{Acknowledgements}
The work of F.F. is in part supported by a startup grant from Melbourne School of Engineering at the University of Melbourne. 

\section*{Author contributions statement}
F.F. is the sole author of the paper.

\section*{Additional information}
The authors declare no competing interests.

\end{document}